\documentclass{amsproc}

\def\C{{\mathbb C}}
\def\N{{\mathbb N}}

\def\R{{\mathbb R}}
\def\Z{{\mathbb Z}}

\def\nv{{{\rm Vol}_{\N}}} 

\def\cd{{\rm codeg}\,}

\newtheorem{theorem}{Theorem}[section]

\theoremstyle{definition}
\newtheorem{dfn}[theorem]{Definition}

\newtheorem{coro}[theorem]{Corollary}
\newtheorem{conj}[theorem]{Conjecture}
\newtheorem{prop}[theorem]{Proposition}
\theoremstyle{remark}
\newtheorem{rem}[theorem]{Remark}

\numberwithin{equation}{section}

\begin{document}

\title{Lattice polytopes with a given 
$h^*$-polynomial}

\author{Victor V. Batyrev}
\address{Department of Mathematics and Physics, 
University of T\"ubingen, Auf der Morgenstelle 10, D-72076 T\"ubingen, Germany}
\email{victor.batyrev@uni-tuebingen.de}





\begin{abstract}
Let $\Delta \subset \R^n$ be an $n$-dimensional lattice polytope. 
It is well-known that 
$h_{\Delta}^*(t) := (1-t)^{n+1} \sum_{k \geq 0} 
|k\Delta \cap \Z^n| t^k $
is a polynomial of degree $d \leq n$ with nonnegative integral  
coefficients. Let $AGL(n, \Z)$ be the group of invertible affine integral 
transformations which naturally acts on $\R^n$. For a given polynomial 
$h^* \in \Z[t]$, we denote by 
$C_{h^*}(n)$ the number $AGL(n, \Z)$-equivalence classes of 
$n$-dimensional lattice polytopes such that $h^* = h_{\Delta}^*(t)$. 
In this paper we show that $\{ C_{h^*}(n) \}_{n \geq 1}$ 
is a monotone increasing sequence which eventually becomes constant.  
This statement follows from a more general combinatorial result 
whose proof uses 
methods of commutative algebra. We give an explict 
description of the sequence   $\{ C_{h^*}(n) \}_{n \geq 1}$ for some 
special polynomials 
$h^*$. 
\end{abstract}

\maketitle

\section{Introduction}

Let $\Delta \subset \R^n$ be an $n$-dimensional lattice polytope, i.e. 
all vertices of $\Delta$ are contained in $\Z^n$. Denote by 
${\rm vol}(\Delta)$ the usual euclidean volume of $\Delta$. 
Define the {\em normalized 
volume} of $\Delta$ to be $$\nv(\Delta):=n! {\rm vol}(\Delta).$$ It is easy 
to see that  $\nv(\Delta)$ is a positive integer.  If $\Delta$ is a 
lattice simplex then $\nv(\Delta) \in \N$  follows from a direct 
calculation (for 
general case  one uses a simplicial subdivision of $\Delta$ into lattice 
simplices). Let us define the {\em codegree} 
of $\Delta$ as  
\[ \cd \Delta := \min \{ k \in \Z_{\geq 0} \; : \; 
|{\rm Int}(k\Delta) \cap \Z^n | \neq 0 \}, \]
where ${\rm Int}(k\Delta)$ denotes the interior of $k\Delta$. 
We notice that  $\cd \Delta \leq  n+1$,  because the sum of vertices 
of any $n$-dimensional 
lattice subsimplex $\Delta' \subset \Delta$ is an element of 
$ {\rm Int}((n+1)\Delta') \cap \Z^n $.     
The nonegative integer  $$\deg \Delta:= n+1 -\cd \Delta $$ will be called 
{\em degree} of $\Delta$. 

Let $\Pi(\Delta) \subset \R^{n+1}$ be  a standard $(n+1)$-dimensional 
pyramid with height $1$ over $\Delta \subset \R^n$, i.e. $\Pi(\Delta)$ is 
 a convex hull of 
the $n$-dimensional polytope $\Delta \cong (\Delta,0) \subset \R^{n+1}$ and 
the lattice vertex  $v=(0, \ldots, 0,1) \in \R^{n+1}$. 

Since ${\rm vol}(\Pi(\Delta)) = \frac{1}{n+1} {\rm vol}(\Delta)$,
we have 
$$\nv( \Pi(\Delta)) = \nv (\Delta).$$ 

By intersecting $k\Pi(\Delta) \subset \R^n$ with hyperplanes $x_{n+1} = m \; 
(1 \leq m  \leq k-1)$ 
we obtain 
\[ |{\rm Int}(k\Pi(\Delta)) \cap \Z^{n+1} | = \sum_{ m =1}^{k-1} 
|{\rm Int}(m\Delta) \cap \Z^n |. \]
The latter implies   $ \cd \Pi(\Delta) =   \cd \Delta + 1$ and 
\[ \deg \Pi(\Delta) =  \deg  
\Delta. \]

\begin{dfn}
{\rm Let $AGL(n, \Z)$ be the group of   affine integral linear transformations
which naturally acts on $\R^n$.  
We denote by $C(V,d,n)$ the number of $AGL(n, \Z)$-equivalence classes 
of $n$-dimensional lattice polytopes such that $\nv \Delta = V$ and 
$\deg \Delta = d$.}  
\end{dfn}

\begin{rem} 
{\rm The fact that there exists only finitely many 
$AGL(n, \Z)$-equiva\-lence classes 
of $n$-dimensional lattice polytopes of fixed volume is  classically 
well-known (see e.g.  \cite{LZ}).} 
\end{rem}

Now we formulate the main result of the paper:

\begin{theorem} 
Let us fix some  integers  $V$ and $d$. Then for 
$$n \geq 4d {2d+ V -1 \choose 2d} $$ every $n$-dimensional 
lattice polytope
$\Delta \subset \R^n$ with  $\nv (\Delta) = V$, $\deg \Delta =d$ 
is a pyramid $\Pi(\Delta')$ over an $(n-1)$-dimensional 
lattice polytope $\Delta' \subset \R^{n-1}$.
\label{main-theo}
\end{theorem} 
\bigskip

\section{Some applications}

Let us  discuss some consequences of Theorem \ref{main-theo}.

\begin{dfn} 
{\rm Let $k$ be a nonnegative integer and $\Delta$ be an $n$-dimensional 
lattice polytope. We define {\em $k$-fold pyramid $\Pi^{(k)}(\Delta)$ over} 
$\Delta$ as 
follows: \\
(i) $\Pi^{(0)}(\Delta) := \Delta$; \\
(ii) $\Pi^{(k+1)} (\Delta) := \Pi( \Pi^{(k)} (\Delta))$. } 
\end{dfn}

\begin{coro} 
For fixed integers $V,d$ there exists a positive constant $\nu=\nu(d,V)$ 
such that 
for $n > \nu$ 
every $n$-dimensional lattice polytope $\Delta$  with $\nv(\Delta) = V$, 
$\deg \Delta = d$ is an  $(n -\nu)$-fold pyramid over a $\nu$-dimensional 
lattice polytope $\Delta'$, i.e. 
\[ \Delta = \Pi^{(n - \nu)} (\Delta'). \]
\end{coro} 

\noindent
\begin{proof} The statement immediately follows  from \ref{main-theo} 
by induction if  we set $$\nu:= 4d  {2d+ V-1 \choose 2d}-1.$$ 
\end{proof} 
\bigskip

Now we show that for fixed integers $V$ and $d$,  one has    
\[ C(V,d,1) \leq   C(V,d,2) \leq  \ldots   \leq  C(V,d,n) \leq  \ldots .\]

\begin{prop} Let $k$ be a nonnegative integer. 
Two $n$-dimensional lattice polytopes $\Delta, \Delta' \subset \R^n$ 
are $AGL(n, \Z)$-equivalent if and only if two 
$k$-fold pyramids 
$\Pi^{(k)}(\Delta), \Pi^{(k)}(\Delta') \subset \R^{n+k}$ are 
$AGL(n+k, \Z)$-equivalent.
\end{prop} 

\noindent
\begin{proof} The case $k =0$ is trivial. 
If  $\Delta$ and $\Delta'$ are $AGL(n, \Z)$-equivalent, then 
the standard embeddings $\Delta \hookrightarrow \Pi^{(k)}(\Delta)$,  
$\Delta' \hookrightarrow \Pi^{(k)}(\Delta')$ and 
 $AGL(n, \Z) \hookrightarrow AGL(n+k, \Z)$ show that 
$\Pi^{(k)}(\Delta)$ and $\Pi^{(k)}(\Delta')$ are  $AGL(n+k, \Z)$-equivalent. 
Now assume that $\Pi^{(k)}(\Delta)$ and $\Pi^{(k)}(\Delta')$ are  
$AGL(n+k, \Z)$-equivalent.  We write 
\[ \Pi^{(k)} (\Delta)= {\rm conv}\{ \Delta, v_1, \ldots, v_k \}, \;\;  
\Pi^{(k)} (\Delta')= {\rm conv}\{ \Delta', v_1', \ldots, v_k' \}. \]
Let $g \in AGL(n+k, \Z)$ be an element such that $g\Pi^{(k)}(\Delta) = 
\Pi^{(k)}(\Delta')$. If $g\Delta = \Delta'$, then we are done. Otherwise
$g\Delta \subset \Pi^{(k)}(\Delta')$ contains  
a vertex $v'_i \in  \Pi^{(k)}(\Delta')$. Therefore, $g\Delta$ is 
itself a pyramid $\Pi(\Theta)$ 
over an $(n-1)$-dimensional face  $\Theta' \subset \Pi^{(k)}(\Delta')$, 
i.e. $\Delta$ is also a pyramid 
$\Pi(\Delta_1)$ over an $(n-1)$-dimensional lattice polytope $\Delta_1 := 
g^{-1}\Theta'$. Analogously $g^{-1} \Delta'$ must contain a vertex 
$v_j  \in \Pi^{(k)}(\Delta)$ and $\Delta'$ is a pyramid $\Pi(\Delta_1')$ 
over an 
 $(n-1)$-dimensional lattice polytope $\Delta_1'$. It remains to 
show that $\Delta_1$ and $\Delta_1'$ are 
$AGL(n-1, \Z)$-equivalent. For this we can apply the same arguments 
as for $\Delta$ and $\Delta'$ and come to either 
$g\Delta_1 = \Delta_1'$, or to another pair of $(n-2)$-dimensional lattice 
polytopes  $\Delta_2$ and $\Delta_2'$ such that $\Delta_1 = 
\Pi(\Delta_2)$ and $\Delta_1' = 
\Pi(\Delta_2')$. After finitely many steps this process terminates. 
\end{proof}

\begin{coro} The correspondence $\Delta \mapsto \Pi^{(k)}(\Delta)$ is injective 
on the set of equivalence classes of lattice polytopes of fixed degree and 
fixed normalized volume considered modulo affine 
integral linear transformations. 
For all natural numbers $V, d, n, k$  one has 
\[ C(V,d,n) \leq C(V,d,n+k), \]
where for sufficiently 
large $n$ the number  $C(V,d,n)$ does not depend on $n$.
\label{mono2}
\end{coro} 
\bigskip

Now consider  two power series 
\[ P(\Delta, t) = \sum_{k \geq 0} |k\Delta \cap \Z^n| t^k, \;\;
Q(\Delta,t) = \sum_{k > 0} |{\rm Int}(k\Delta) \cap \Z^n |t^k. \]
It is well-known that $P(\Delta,t)$ and $Q(\Delta,t)$ are rational 
functions such that  
\[ Q(\Delta,t) = (-1)^{n+1} P(\Delta,t^{-1}), \;\; 
 P(\Delta,t) = \frac{ h_{\Delta}^*(t)}{(1-t)^{n+1}},\]
where $ h_{\Delta}^*(t)$ is  a  polynomial with 
nonnegative integral coefficients  satisfying the conditions 
\[ h_{\Delta}^*(0)=1, \;\; 
 h^*_{\Delta}(1)= \nv(\Delta). \]

\begin{dfn} 
{\rm We call the polynomial 
\[ h^*_{\Delta}(t) := (1 -t)^{n+1} P(\Delta,t) \]
$h^*$-{\em polynomial} of the lattice polytope $\Delta$ 
(this definition is inspired by the notion of $h^*$-vector 
considered by Stanley in \cite{Sta93}). }
\end{dfn} 

\begin{rem} 
{\rm 
It follows from the formula  
$$Q(\Delta,t) = \frac{t^{n+1} h^*_{\Delta}(t^{-1})}{(1-t)^{n+1}}$$ 
that $\cd \Delta$ equals 
the multiplicity of the root $t=0$ of the polynomial 
$t^{d+1}h^*_{\Delta}(t^{-1})$ and 
\[ \deg \Delta = \deg h_{\Delta}^*(t). \] 
Since
 \[ P(\Pi(\Delta), t) = \frac{P(\Delta, t)}{1-t},  \]
we have also
\[ h_{\Delta}^*(t) = h_{\Pi(\Delta)}^*(t). \]
} 
\end{rem}

\begin{dfn} 
{\rm 
Let $h^*(t) \in \Z[t]$ be a polynomial with nonegative integral coefficients.
We denote by  $C_{h^*}(n)$ the number of $AGL(n, \Z)$-equivalence classes 
of $n$-dimensional lattice polytopes $\Delta$  such that 
$h^*_{\Delta} = h^*$.} 
\end{dfn}

Since there exist only finitely many polynomials $h^*(t)$ 
with nonnegative integral coefficients of given degree $d$ and given sum 
of coefficients 
$h^*(1) = V$ we obtain 
\[ C(V, d,n) = \sum_{h^*(1) = V, {\rm deg}\, h^* =d} C_{h^*}(n). \] 
Theorem \ref{main-theo} together with \ref{mono2} 
implies the following:  

\begin{coro} 
Let $h^*(t) \in \Z[t]$ be a polynomial of degree $d$. Then for sufficiently 
large $n$ every $n$-dimensional lattice polytope $\Delta \subset 
\R^n$ such that 
$h_{\Delta}^* = h^*$ is a pyramid $\Pi(\Delta')$ 
over a $(n-1)$-dimensional lattice polytope $\Delta' \subset \R^{n-1}$, i.e.  
one obtains a monotone sequence 
\[  C_{h^*}(1) \leq  C_{h^*}(2) \leq \cdots \leq  C_{h^*}(n) \leq \cdots \]
which eventually becomes constant. 
\end{coro} 
\bigskip

\section{The proof}

The proof of Theorem \ref{main-theo} uses methods of commutative algebra 
in the spirit of \cite{MS}. 

Let $\sigma_{\Delta}$ be the  $(n+1)$-dimensional cone in $\R^{n+1}$ 
over $\Delta$: 
\[ \sigma_{\Delta} := \{ (r\Delta, r) \in \R^{n+1} \; : \; r \in \R_{\geq 0} \} 
\subset \R^{n+1} . \] 
The set $M_{\Delta}:= \sigma_{\Delta}\cap \Z^{n+1}$ of all lattice points 
in the cone  $\sigma_{\Delta}$ is a monoid with respect to sum. 
Moreover,  $M_{\Delta}$ is a {\em graded} monoid with respect to the 
$(n+1)$-th coordinate, i.e. the degree of a lattice 
point $(m,k) \in \Z^n \times \Z$ equals  $k$. 
We define $S_{\Delta}:= \C[M_{\Delta}]$ to be the graded semigroup $\C$-algebra
of the graded monoid $M_{\Delta}$. Since the $k$-th homogeneous component  
$M_{\Delta}^k$ of  $M_{\Delta}$ has form  
\[ M_{\Delta}^k = \{ (m,k) \in \Z^n \times \Z \; : \; m \in k \Delta \} \]
we have $$\dim_{\C}   S_{\Delta}^k = | M_{\Delta}^k| =  
|k\Delta \cap \Z^n|.$$

This allows to interpret the power series  $P(\Delta, t)$ as  a 
Hilbert-Poincar{\'e} series of the graded commutative $\C$-algebra
$S_{\Delta} = \bigoplus_{k \geq 0} S_{\Delta}^k$, i.e.   
\[ P(\Delta, t) = \sum_{k \geq 0} (\dim_{\C}  S_{\Delta}^k) t^k. \]
We remark that $M_{\Delta}$ is a finitely generated 
graded monoid, i.e. $S_{\Delta}$ is a graded finitely generated 
$\C$-algebra.

\begin{dfn} 
{\rm We call a set  ${\mathcal X} := 
\{x_1, \ldots, x_l \} \subset \sigma_{\Delta} \cap \Z^{n+1}$ 
a {\em minimal generating set  
of  the monoid}  $M_{\Delta}$ if  ${\mathcal X}$ generates the 
monoid  $M_{\Delta}$ and every lattice point $x_i \in {\mathcal X}$ cannot
be represented as a linear combination of  
${\mathcal X} \setminus \{ x_i \}$ with nonnegative integral coefficients.} 
\end{dfn}

\begin{rem} {\rm We note that every minimal generating set  
${\mathcal X} := 
\{x_1, \ldots, x_l \}$ of $M_{\Delta}$ must contain  
the set $M_{\Delta}^1$ of all lattice points in $\sigma_{\Delta}$ 
of degree $1$, because these lattice points cannot be represented 
as nonnegative integral linear combination of other lattice points 
in $M_{\Delta}$. So we have  $M_{\Delta}^1 \subset 
{\mathcal X}$ and $l \geq |M_{\Delta}^1| = 
|\Delta \cap \Z^n | \geq n+1$.
} 
\label{deg1}
\end{rem}

\begin{dfn} 
{\rm 
We associate with each lattice point $x_i \in {\mathcal X}$ a variable 
$X_i$ and denote by $A$ the polynomial algebra 
\[  A:=  \C[X_1, \ldots, X_l ]. \] 
The grading of  $A$ is defined  by the grading of the lattice 
points in  ${\mathcal X}$:
\[ \deg X_i := \deg x_i , \;\; i =1, \ldots, l. \]
} 
\end{dfn} 

\begin{rem} {\rm 
Now the finitely generated $\C$-algebra $S_{\Delta}$ can be written as  
\[ S_{\Delta} \cong  A/I \]
where  $I$ is the homogeneous ideal in  
$A$ generated by binomials
\[B= B(R):=  X_{i_1}^{a_1} \ldots  X_{i_s}^{a_s} -   X_{j_1}^{b_1} \ldots  
X_{j_r}^{b_r}, \]
corresponding to linear relations
\[ R\;\; : \;\; a_1 x_{i_1} + \cdots + a_s x_{i_s} = b_1 x_{j_1} + \cdots + 
b_r x_{j_r},   \]
\[ \{ i_1, \ldots, i_s \} \cap \{ j_1, \ldots, j_r \} = \emptyset, 
\;\; a_i, b_j \in \Z_{>0}. \]
} 
\end{rem}

The key  observation in the proof of Theorem \ref{main-theo} is 
the following statement:

\begin{prop} 
Let ${\mathcal X} := 
\{x_1, \ldots, x_l \}$ be  a minimal generating set of $M_{\Delta}$ and 
let $\{ B_1, \ldots,B_p \}$ be a set of binomials generating 
the ideal  $I$. 
Then $\Delta = \Pi(\Delta')$ for some $(n-1)$-dimensional lattice 
polytope $\Delta'$ if and only if there exists a lattice point 
$x_i \in M_{\Delta}^1$ such that the corresponding variable $X_i$ 
does not appear in any of binomials $B_1, \ldots, B_p$. In the latter
case $x_i$ is the vertex of the pyramid  $\Pi(\Delta')$ and  
$S_{\Delta} \cong S_{\Delta'}[X_i]$ is a polynomial ring over $S_{\Delta'}$.
\label{poly}
\end{prop}

\noindent 
\begin{proof}  ``$\Rightarrow$'': 
Let  $\Delta = \Pi(\Delta') \subset \R^n$ for some 
$(n-1)$-dimensional lattice polytope $\Delta'$. We can assume that 
all points of the  polytope $\Delta'$ have zero last $n$-th coordinate and 
the vertex $v$ 
of the pyramid $\Pi(\Delta')$ is the lattice point $(0, \ldots, 0,1) 
\in \R^n$. In this case, all lattice points in the monoid $M_{\Delta} 
\subset \Z^{n+1}$ have nonnegative $n$-th coordinates. 
Since $(v,1) \in M_{\Delta}^1$, by \ref{deg1}, we obtain 
 $(v,1) \in {\mathcal X}$, i.e., $(v,1) = x_i$ 
for some $1 \leq i \leq l$. Assume that for some $j \in \{1, \ldots, p \}$ 
the  variable 
$X_i$ corresponding to $x_i$ appears in the binomial
\[ B_j =  X_{i_1}^{a_1} \ldots  X_{i_s}^{a_s} -   X_{j_1}^{b_1} \ldots  
X_{j_r}^{b_r}. \]
Without loss of generality we can assume  that 
$i = i_1$, $(a_1 > 0)$. Then one has an integral linear relation 
 \[  a_1 x_i + a_2 x_{i_2} + \cdots + a_s x_{i_s} = b_1 x_{j_1} + \cdots + 
b_r x_{j_r}, \]
\[ \{ i, i_2, \ldots, i_s \} \cap \{ j_1, \ldots, j_r \} = \emptyset. \]
Since $x_i =(v,1) \in \Z^{n+1}$ has positive $n$-th coordinate 
and all lattice points in  ${\mathcal X}$ have nonnegative $n$-th 
coordinate, we obtain that there exists $q \in \{j_1, \ldots, j_r \}$ 
such that $x_q$ also has a positive $n$-th coordinate. Using the splitting
\[ M_{\Delta} = M_{\Delta'} \oplus \Z x_i ,\]
we can write $x_q = x_q' + cx_i$ where $x_q' \in  M_{\Delta'}$ and $c>0$
(all elements of  $M_{\Delta'}$ have zero $n$-th coordinate). 
We remark that  $x_q' \neq 0$ because otherwise we would have a contradiction
to minimality of ${\mathcal X}$ or to the condition  
$\{ i, i_2, \ldots, i_s \} \cap \{ j_1, \ldots, j_r \} = \emptyset$.   
Therefore  $x_q'$ is a positive linear combination of some lattice points 
from ${\mathcal X}_0:= {\mathcal X} \cap M_{\Delta'}$. This also contradicts 
the minimality of ${\mathcal X}$. Thus $x_i$ does not appear in any of 
binomials $B_1, \ldots, B_p$. 

``$\Leftarrow$'': Let $x_i \in  {\mathcal X}$ be a lattice point such that 
the corresponding variable $X_i$ does not appear in any of binomials
$B_1, \ldots, B_p$. Let $M_{\Delta}'$ be the monoid generated 
by  ${\mathcal X}':= 
{\mathcal X} \setminus \{ x_i \}$. Then $S_{\Delta} \cong A/(B_1, \ldots, B_p)$ 
is isomorphic 
to the polynomial ring $\C[M_{\Delta}'][X_i]$. In particular, we have 
\[ {\rm Krull}\, \dim \C[M_{\Delta}' ] =  {\rm Krull}\,\dim S_{\Delta} - 1 
= n. \] 
Therefore all lattice points from  ${\mathcal X}'$ belong to a $n$-dimensional 
linear subspace $L \subset \R^{n+1}$. The subspace $L$ 
cuts the cone $\sigma_{\Delta}$ 
along its $n$-dimensional face $\Theta$, because all generators 
${\mathcal X}$ of $M_{\Delta}$ except $x_i$ belong 
to $L$.
Denote by $H$ the affine hyperplane in $\R^{n+1}$ consisting of all 
points $y \in \R^{n+1}$ whose  $(n+1)$-th coordinate equals $1$. 
Then $x_i \in H$, $\Delta = \sigma_{\Delta} \cap H$ and 
$\Delta' := \Theta \cap H$ is a convex polytope such that $\Delta = 
{\rm conv} ( x_i, \Delta')$. Therefore all vertices of $\Delta'$ are 
lattice points and $M_{\Delta'} = M_{\Delta} \cap L$. Since there is 
no nontrivial relations between $x_i$ and  $M_{\Delta'}$ we have 
\[  M_{\Delta} =   M_{\Delta'} \oplus \Z x_i. \]
Thus $\Delta$ is a standard pyramid over $\Delta'$ with vertex $x_i$. 
\end{proof}
\bigskip

By well-known result of Hochster \cite{Ho}, 
$S_{\Delta}$ is a Cohen-Macaulay ring of Krull dimension $n+1$. 
We consider a minimal  graded free resolution of $S_{\Delta}$ as an
$A$-module: 

\[ P^{\bullet} \;\; : \;\; 
0 \to P_{l -n-1} \to \cdots \to P_1 \to P_0 \to S_{\Delta} \to 0 \]
where $P_i$ is a free graded $A$-module of rank $\alpha_i$. The minimality
of the resolution means that all elements of the $\alpha_{i-1} \times 
\alpha_{i}$-matrix of the differential $d_i\; : \; P_i \to P_{i-1}$ 
are contained in the maximal homogeneous ideal $\mu$ of $A$. It is well-known 
that such a free resolution is uniquely determined up to isomorphism. 
Moreover, one has 
\[ \alpha_i = \dim_{\C} {\rm Tor}_i^A(\C, S_{\Delta}) \]
where $\C$ is considered as $A$-module via the isomorphism $\C \cong A/\mu$. 
In particular, we have  $\alpha_0 =1$ and $\alpha_1$ equals the minimal number
of generators of the ideal $I$, i.e. $\alpha_1 = \dim_{\C} I/\mu I$. 

Let $\{ y_0, y_1, \ldots, y_n \} \subset S_{\Delta}^1$be 
a maximal $S_{\Delta}$-regular sequence consisting of elements
of degree $1$. The existence of such a regular sequence follows from the fact 
that $ S_{\Delta}^1$ generates a primary ideal in $S_{\Delta}$ whose 
radical is the maximal homogeneous ideal  in $S_{\Delta}$. 
Using the surjective homomorphism $A \to S_{\Delta}$, 
we can choose  $n+1$ linearly  independent 
homogeneous linear
forms $Y_0, Y_1, \ldots, Y_n \in A= \C[X_1, \ldots, X_l]$ whose images 
in $S_{\Delta}$ coincide with $y_0, y_1, \ldots, y_n$. 
Let us  define 
\[ R:=  S_{\Delta}/ (y_0, y_1, \ldots, y_n),\;\; 
\overline{I} = I/(Y_0, Y_1, \ldots, Y_n)I,  \]
\[ \overline{A}:= A/(Y_0, Y_1, \ldots, Y_n), \;  \overline{\mu}:= 
\mu/(Y_0, Y_1, \ldots, Y_n), \]
\[ \overline{P_i} = P_i/(Y_0, Y_1, \ldots, Y_n)P_i, \;\; i =0, \ldots, l-n-1. \]

\begin{prop} 
The complex 
$\overline{P^{\bullet}} = P^{\bullet} \otimes_A \overline{A}$ 
is exact so that we can consider
\[ \overline{P^{\bullet}}\;: \; 0 \to    \overline{P_{l - n -1}} 
\to  \cdots \to  \overline{P_{1}}  \to \overline{P_{0}} \to R \to 0\]
as a minimal graded free resolution of $R \cong \overline{A}/ 
\overline{I}$ as a  $\overline{A}$-module.
In particular, the numbers and the degrees of the minimal 
generators of the ideals $I \subset A$ and $\overline{I} 
\subset \overline{A}$ are the same.
\label{2same}
\end{prop} 

\noindent
\begin{proof} 
The $i$-th cohomology of  the complex 
$\overline{P^{\bullet}} = P^{\bullet} \otimes_A \overline{A}$ equals 
\[ {\rm Tor}_{i}^A ( \overline{A}, S_{\Delta}). \]
On the other hand, one can compute 
${\rm Tor}_{i}^A ( \overline{A}, S_{\Delta})$
using the  Koszul complex on the  elements $Y_0, Y_1, \ldots, 
Y_n$ which is a graded free resolution of  $\overline{A}$ over 
$A$. Since the sequence $Y_0, Y_1, \ldots, Y_n$ is 
$S_{\Delta}$-regular ${\rm Tor}_{i}^A ( \overline{A}, S_{\Delta}) =0$ 
for all $i >0$ and  ${\rm Tor}_{0}^A ( \overline{A}, S_{\Delta}) =R$.  
Therefore $\overline{P^{\bullet}} $ is exact. The minimality of the 
graded free resolution $\overline{P^{\bullet}} $ follows from the fact 
that the maximal homogeneous ideal $\overline{\mu} \subset \overline{A}$ 
is the image under the surjection $A \to \overline{A}$ 
of the  maximal homogeneous ideal ${\mu} \subset {A}$. 
Since 
\[ P_1/\mu P_1 \cong \overline{P_1}/\overline{\mu} \overline{P_1} \]
we obtain that  the numbers and the degrees of the minimal 
generators of two ideals $I \subset A$ and $\overline{I} 
\subset \overline{A}$ are the same.
\end{proof}

\begin{prop} 
Let $\Delta$ be an $n$-dimensional lattice polytope such that  $V= 
\nv(\Delta)$. Then the monoid $M_{\Delta}$ contains a minimal
generating subset ${\mathcal X}$ containing $l \leq V + n$ elements which all 
have  degree $\leq d$.
\label{ineq}
\end{prop}

\noindent
\begin{proof} We remark that the coefficients $h^*_i$ of the 
$h^*$-polynomial $h^*_{\Delta}$ are equal to the dimensions 
of the homogeneous 
components of the graded artinian ring 
\[ R = R^0 \oplus R^1 \oplus \cdots \oplus R^d, \]
i.e. $h^i = \dim_{\C} R^i$ and $h^*_{\Delta}(1) = V$.
Since  $S_{\Delta}$ is Cohen-Macaulay we obtain that  
the $\C$-algebra $S_{\Delta}$ is a free module of rank 
$V$ over the polynomial ring $B:= \C[y_0,y_1, \ldots, y_n]$
( $S_{\Delta}$ is an integral extension of  $\C[y_0,y_1, \ldots, y_n]$). 
Consider ${\mathcal Z} := M_{\Delta}^1 \cup \{ z_1, \ldots, z_m \}$ where
$z_1, \ldots, z_m$ are some lattice points of degree $k$ $(2 \leq k \leq d)$ 
in $M_{\Delta}$ such that their images $\overline{z_1} \ldots, \overline{z_m}$ 
in $R$ form a basis of the $\C$-vector space $\bigoplus_{i=2}^d R^i$.
We note that  ${\mathcal Z}$ generates the $\C$-algebra $S_{\Delta}$ (and hence
also the monoid $M_{\Delta}$), because the $\C$-algebra generated 
by ${\mathcal Z}$ contains the polynomial ring  
$B$ and all generators of $S_{\Delta}$ as a finite $B$-module. 
Thus we have 
\[ | {\mathcal Z}| = | M_{\Delta}^1| + \sum_{i \geq 2} h^*_i. \]
Since 
$h_0^* =1, \;\; h^*_1 =  | M_{\Delta}^1| - n -1$
the generating subset  ${\mathcal Z}$ contains
$n + \sum_{i \geq 0} h^*_i = n + V$ 
lattice points and all elements in  ${\mathcal Z}$ have degree $\leq d$. 
\end{proof}

\begin{prop} 
All binomials in a minimal generating set for the ideal $I \subset A$ 
have degree 
at most $2d$. 
\label{bi-deg}
\end{prop}

\noindent
\begin{proof} By \ref{2same}, it is sufficient to prove the same 
statement for the ideal $\overline{I} 
\subset \overline{A}$. Let  $ \overline{A}^k$ and $\overline{I}^k$ denote 
$k$-th  homogeneous components of the graded ring  $\overline{A}$ and 
its homogeneous ideal $\overline{I}$. Since for  $k \geq  d+1$
one has  $\overline{A}^k =\overline{I}^k,$ 
it is sufficient to show that for all $i >0$ the homogeneous component 
 $\overline{A}^{2d+i}$ is generated by products 
\[  \overline{A}^{d -j} \cdot \overline{A}^{d+i +j }  \;\; (0 \leq j 
\leq d -1). \]
The latter follows from the fact that $A$ (and hence also $\overline{A}$) 
is generated by elements of degree $\leq d$ (see \ref{ineq}). 
\end{proof} 

\begin{prop} 
The number of binomials  in a minimal generating set for the ideal 
$I \subset A$ 
is not greater than 
\[ { 2d +V-1 \choose 2d }. \]
\label{bi-num} 
\end{prop} 

\noindent
\begin{proof}  By \ref{2same}, it is sufficient to prove the same 
statement for a minimal generating set of  the ideal $\overline{I} 
\subset \overline{A}$. By \ref{bi-deg}, the number of 
minimal generators of 
$\overline{I} \subset \overline{A}$ is not greater 
than the dimension of the space of all polynomials in $l-n-1$ variables
of degree $\leq 2d$. The latter is not 
 greater than 
\[ { 2d + l-n-1 \choose 2d }, \]
because the maximum of this dimension is attained if all 
$l-n-1$ variables in $\overline{A}$ have degree $1$. It remains
to apply the inequality $l \leq V + n$ (see \ref{ineq}). 
\end{proof}
\bigskip
\bigskip

{\sc Proof of Theorem \ref{main-theo}.}
Let $B_1, \ldots, B_p$ be a minimal generating set of binomials
for the ideal $I$. 
The number of different variables from $\{ X_1, \ldots, X_l \}$ 
appearing in a binomial $B_i$ is obviously not greater
than $2 \deg B_i$ 

By \ref{bi-num}, we have $p \leq  { 2d +V-1 \choose 2d }$. 
By \ref{bi-deg}, we have $\deg B_i \leq 2d$. 
Therefore the number of variables from 
$\{ X_1, \ldots, X_l \}$ appearing in at least 
one binomial relation is not greater than 
\[ 4d { 2d +V-1 \choose 2d } .\]
If $$ n \geq  4d { 2d +V-1 \choose 2d }, $$
then by \ref{deg1} 
\[  |\Delta \cap \Z^n | >   
4d { 2d +V -1 \choose 2d }
\]
and hence there exists a lattice point $x_i \in  \Delta \cap \Z^n$ such that 
the corresponding variable $X_i$ does not appear in any 
binomial relation $B_1, \ldots, B_p$. 
Hence, by \ref{poly}, $\Delta = \Pi(\Delta')$ for some
$(n-1)$-dimensional lattice polytope. 

\hfill $\Box$ 
\bigskip

\section{Some examples and conjectures}

We remark that the estimate for $n$ in Theorem \ref{main-theo} 
is far from being optimal.

Using the complete classification of lattice polytopes with linear 
$h^*$-polynomial \cite{BN}, one can immediately  explicitly compute 
the numbers $C(V,1,n) = C_{h^*}(n)$ where 
$h^* = h_{\Delta}^* = 1 + (V-1)t $.

\begin{prop} 
Let $p(n,V)$ be 
 the number of integral
solutions of  the equation
\[ \sum_{i=1}^n k_i = V, \;\; 
0 \leq k_1 \leq \cdots \leq k_n. \]
Then  $C(4,1,1) = 1$, $C(4,1,n) = p(n,4) +1$ (if   $n \geq 2)$, and  
$C(V,1,n) = p(n,V)$, (if $V \neq 4, n \geq 1)$.  
In particular, the  monotone sequence $C(V,1,n)$ becomes  constant for 
$n \geq V$.
\end{prop} 

\noindent
\begin{proof} By  \cite{BN}, every $n$-dimensional 
lattice polytope with linear $h^*$-poly\-no\-mial is either a $(n-2)$-fold 
pyramid of the lattice triangle $T \subset \R^2$ with vertices 
$(0,0), (2,0), (0,2)$, or a Lawrence prism with heights $k_1, \ldots, 
k_n$ $( \sum_{i =1}^b k_i = V)$. Up to an $AGL(n, \Z)$-isomorphism, 
we can always assume that $0 \leq k_1 \leq \cdots \leq k_n$. 
\end{proof}

It is interesting to understand in general which polynomials in 
$\Z[t]$ can be 
realised as $h^*$-polynomials of lattice polytopes. In this connection we 
propose the following conjecture:

\begin{conj} 
Let $h^*= \sum_{0\leq i \leq d} h^*_i t^i \in \Z[t]$ be 
a polynomial of degree $d$ with nonnegative 
coefficients. 
If $h^* = h_{\Delta}^*$ for some $n$-dimensional lattice polytope $\Delta$, 
then $h^*(1)= \nv \Delta$ is bounded by some constant 
depending only on the leading coefficient $h^*_d$ of $h^*$.
\label{conj1}     
\end{conj} 

By a theorem 
of Hensley \cite{H},  \ref{conj1} is known to be true for $n =d$. 
Obviously \ref{conj1} is true for arbitrary $n$ if $h^*$ is 
a linear polynomial. 
For quadratic polynomials $h^*$,  we expect the following more precise
statement:

\begin{conj} 
For any $n$-dimensional lattice polytope with quadratic 
$h^*$-polynomial one has:  
\[ h_1^* \leq 3 h_2^* + 4. \]
\label{conj2}
\end{conj} 

For $n =2$  this conjecture is known to be true by a theorem 
of Scott \cite{Sc76} (see also \cite{HS}).

\bibliographystyle{amsalpha}

\begin{thebibliography}{A}

\bibitem [BN]{BN} V. Batyrev, B. Nill, {\em Multiples of lattice 
polytopes without interior lattice points}, Preprint, math.CO/0602336.


\bibitem [HS]{HS}  C. Haase, J. Schicho: 
{\em Lattice polygons and the number 2i+7}, 
math.CO/0406224

\bibitem[He]{H} D. Hensley, {\em 
Lattice vertex polytopes with interior lattice points}, 
Pacific J. Math. 105 (1983), no. 1, 183--191.

\bibitem[Ho]{Ho} M. Hochster,{\em Ring of invariants of tori, 
Cohen-Macaulay rings generated by monomials, and polytopes}, 
Ann. of Math. (2){\bf 96} (1972), 328-337.


\bibitem[MS]{MS} E. Miller, B. Sturmfels, {\em Combinatorial 
commutative algebra}, Springer 2005. 


\bibitem[Sc]{Sc76} R.P.Scott, {\em On convex lattice polygons}, 
Bull. Austral. Math. Soc. 15(1976),
p.395-399

\bibitem[St]{Sta93} R.P. Stanley, {\em A monotonicity property of 
$h$-vectors and $h^*$-vectors}. Eur. J. Comb. \textbf{14}, 251-258 (1993)

\bibitem[LZ]{LZ} J.S. Lagarias, G.M. Ziegler, {\em Bounds for lattice 
polytopes containing a fixed number of interior 
lattice points in a sublattice }, Can. J. Math. {\bf 43}, 1022-1035
(1991).  

\end{thebibliography}

\end{document}